\input amstex
\documentstyle{amams} 
\document
\annalsline{153}{2001}
\received{October 2, 1996}
 \startingpage{1}
\font\tenrm=cmr10
\def\joinrel{\mathrel{\mkern-4mu}}
\def\relbar{\mathrel{\smash-}}
\def\lrar{\relbar\joinrel\relbar\joinrel\relbar\joinrel\relbar\joinrel\rightarrow}
\def\bra#1{\buildrel{#1}\over {\longrightarrow}}
\catcode`\@=11
\font\twelvemsb=msbm10 scaled 1100

\font\ninemsb=msbm10 scaled 800
\newfam\msbfam
\textfont\msbfam=\twelvemsb  \scriptfont\msbfam=\ninemsb
  \scriptscriptfont\msbfam=\ninemsb
\def\msb@{\hexnumber@\msbfam}
\def\Bbb{\relax\ifmmode\let\next\Bbb@\else
 \def\next{\errmessage{Use \string\Bbb\space only in math
mode}}\fi\next}
\def\Bbb@#1{{\Bbb@@{#1}}}
\def\Bbb@@#1{\fam\msbfam#1}
\catcode`\@=12

 \catcode`\@=11
\font\twelveeuf=eufm10 scaled 1100
\font\teneuf=eufm10
\font\nineeuf=eufm7 scaled 1100
\newfam\euffam
\textfont\euffam=\twelveeuf  \scriptfont\euffam=\teneuf
  \scriptscriptfont\euffam=\nineeuf
\def\euf@{\hexnumber@\euffam}
\def\frak{\relax\ifmmode\let\next\frak@\else
 \def\next{\errmessage{Use \string\frak\space only in math
mode}}\fi\next}
\def\frak@#1{{\frak@@{#1}}}
\def\frak@@#1{\fam\euffam#1}
\catcode`\@=12

\font\fr=eufm10

\define\a{\alpha}
\predefine\barunder{\b}
\redefine\b{\beta}

\predefine\dotunder{\d}
\redefine\d{\delta}
\define\e{\varepsilon}
\define\equdef{\overset\text{def}\to=}

\define\g{\gamma}

\define\lb\{{\left\{}
\define\la{\lambda}

\define\lm{\limits}
\define\lra{\longrightarrow}

\define\n{\nabla}

\define\oper{\operatorname}

\define\ov{\overline}

\define\rb\}{\right\}}

\define\vp{\varphi}
\define\wh{\widehat}
\define\wt{\widetilde}
\define\x{\times}

\define\({\left(}
\define\){\right)}
\define\[{\left[}
\define\]{\right]}
\define\<{\left<}
\define\>{\right>}
\def\slantline#1#2#3#4#5{\hbox to 0pt{
}}


\def\SD{{\Cal D}}
\def\SE{{\Cal E}}

\def\SI{{\Cal I}}

\def\SM{{\Cal M}}

\def\SS{{\Cal S}}
\def\ST{{\Cal T}}
\def\SU{{\Cal U}}


\def\BI{{\bold I}}

\def\BK{{\bold K}}

\def\BM{{\bold M}}

\def\BP{{\bold P}}

\def\BT{{\bold T}}

\def\BZ{{\bold Z}}


\def\bbc{{\Bbb C}}

\def\bbp{{\Bbb P}}

\def\bbr{{\Bbb R}}

\def\bbz{{\Bbb Z}}


\def\equdef{\overset{\text{def}}\to{=}}

\def\gg{{\hbox{\fr g}}}

\def\tv{\widetilde V}

\def\tfrac#1#2{{\tsize{\frac{#1}{#2}}}}
\def\ts{\BT}
\def\ps{\BP}
\def\ov{\overline}
\def\cct{{\Cal  T} \!\!\! l}

\def\bbc{{\Bbb C}}

\def\bbp{{\Bbb P}}

\def\bbr{{\Bbb R}}

\def\bbz{{\Bbb Z}}
\def\arr{\longrightarrow}
\def\ct{{\Cal T}}
\def\cs{{\Cal S}}
\def\co{{\Cal O}}

\def\ce{{\Cal E}}

\def\cz{{\Cal Z}}
\def\Hom{{\text{Hom}}}

\def\g{\gamma}
\def\pr{\text{pr}}

\def\dim{{\rm dim}}

\def\D{\Delta}
\def\l{\lambda}
\def\n{\nu}

\def\vf{\varphi}
\def\b{\beta}
\def\pr{\text{pr}}

\def\co{{\Cal O}}

\def\O{\overline}
\define\bd{\partial}
\define\Kahler{K{\"a}hler\ }

\define\Lemma#1{\medskip\noindent{\bf Lemma {#1.\  }}}
\define\Theorem#1{\medskip\noindent{\bf Theorem #1.\  }}
\define\Prop#1{\medskip\noindent{\bf Proposition #1.\  }}
\define\Cor#1{\medskip\noindent{\bf Corollary #1.\  }}
\define\Note#1{\medskip\noindent{\bf Note #1.\  }}
\define\Def#1{\medskip\noindent{\bf Definition #1.\  }}
\define\Remark#1{\medskip\noindent{\bf Remark #1.\  }}
\define\Ex#1{\medskip\noindent{\bf Example #1.\  }}
\define\Arr#1{\overset{#1}\to{\longrightarrow}}

\def\Ci{C^{\infty}}

\def\1#1{\varphi_#1}
\def\0#1#2{\varphi^#1_#2}
\def\2#1{\Phi_#1}
\def\3#1{\ST_#1}

\def\piy{\pi^*_Y}
\def\pix{\pi^*_X}
\def\pt{\BP_t}
\def\bp#1{\BP_#1}
\def\bt{\BT_t}

\def\xx{X\x X}
\def\yx{Y\x X}
\def\graph{\oper{graph}}

\def\pr{\oper{pr}}
\def\crr{{\rm Cr}}
\def\spt{\oper{spt}}

\def\cur{\SD'}
\def\Stab{\SS}
\def\Ustab{\SU}

\title{Finite volume flows and Morse theory} 
\acknowledgements{The research of both authors was partially supported by the NSF.}
 \twoauthors{F. Reese Harvey}{H. Blaine Lawson, Jr.}
   \institutions{Rice University, Houston, TX\\
{\eightpoint {\it E-mail address\/}:  harvey\@math.rice.edu}\\
\vglue6pt
State University of New York at Stony Brook, Stony Brook, NY\\
{\eightpoint {\it E-mail address\/}: blaine\@math.sunysb.edu}}

\advance\sectioncount by -1
\section{Introduction}
  In this paper we present a new
 approach to Morse theory based on the de~Rham-Federer theory of currents. The full classical theory is
derived in a transparent way.   The methods carry over uniformly to the equivariant
and the holomorphic settings.  Moreover, the methods are substantially
stronger than the classical ones and have interesting applications to
geometry.  They lead, for example, to formulas relating characteristic
forms and singularities of bundle maps.

The ideas came from addressing the following.

\demo{Question} Given a smooth flow $\vf_t:X\to X$ on a manifold $X$,
when does the limit $$\BP(\a) \equiv \lim_{t\to\infty}\vf_t^*\a$$
exist for a given smooth differential form $\a$?
\enddemo

We do not demand that this limit be smooth.  As an example consider the
gradient flow of a linear function restricted to the unit sphere $S^n\subset
\bbr^{n+1}$.  For any $n$-form $\a$, $\BP(\a)= c [p]$ where
$c=\int_{S^n}\a$ and $[p]$ is the point measure of the bottom critical
point $p$.
\footnote{**}{When $\a$ has compact support in ${\scriptstyle\bbr}^n=S^n-\{\infty\}$ this
simple example is the useful ``approximate identity'' argument widespread
in analysis.}

Our key to the problem is the concept of a {\it finite-volume flow} --
a flow for which the graph $T$  of the relation $x \prec y$,
defined by the forward motion  of the flow, has finite  volume in $\xx$.
Such flows are abundant.
They include generic gradient flows, flows for which the graph is analytic
(e.g., algebraic flows),
and any flow with fixed points on $S^1$.

 We show that
for any flow of finite volume, the limit $\BP$ exists for all $\a$   and
defines a  continuous linear operator
$$
\BP : \ce^*(X) \ \arr\ {\cur}^*(X)
$$
from smooth forms to generalized forms, i.e., currents.  Furthermore, this operator is
chain homotopic to the  inclusion $\BI: \ce^*(X) \hookrightarrow
{\cur}^*(X)$.  More precisely, using the graph $T$, we shall construct a
continuous operator $ \BT : \ce^*(X) \ \arr\ {\cur}^*(X) $
 of degree $-1$ such that
$$
d \circ \BT + \BT \circ d \ =\ \BI - \BP.
\tag{FME}
$$
By   de Rham [deR], $\BI$ induces an  isomorphism in
cohomology.  Hence so does~$\BP$.

Now let   $f:X \to \bbr$ be a Morse function
 with   critical set ${\rm Cr}(f)$, and suppose there
is a Riemannian metric on $X$ for which the gradient flow
$\vp_t$ of $f$ has the following properties:

\medbreak
\item{(i)}    $\vp_t$ is of finite volume.
\smallbreak
\item{(ii)}  The stable and unstable manifolds, $S_p$, $U_p$ for
$p\in {\rm Cr}(f)$, are of finite volume in $X$.
\smallbreak
\item{(iii)}  $p\prec q \ \Rightarrow \ \lambda_p < \lambda_q$ for all
$p,q \in {\rm Cr}(f)$,  where $\lambda_p$ denotes the index of $p$
and where $p\prec q $ means  there is a piecewise flow line
connecting $p$ in forward time to $q$.
\medbreak
\noindent
Metrics yielding such gradient flows always exist.  In fact  Morse-Smale
gradient systems have these properties [HL$_4$]. Under the  hypotheses
(1)--(3) we prove that the operator $\BP$ has the following simple form
$$
\BP(\a) \ = \ \sum_{p\in {\rm Cr}(f)}{\tsize\left( \int_{U_p} \a
\right)}[S_p]
$$
where
$
\int_{U_p} \a     = 0 $ if $\deg \a \ne \lambda_p.
$
Thus $\BP$ gives a retraction
$$
\BP : \ce^*(X) \ \arr\
\Stab_f \ \equdef\ \text{span}_\bbr \biggl \{ [S_p]\biggr\}_{p \in {\rm Cr}(f)}
$$
onto the finite-dimensional subspace of currents spanned by the stable
manifolds of the flow.  Using (FME) we see that $\Stab_f$ is
$d$-invariant and  show that
$
H^*\left(\Stab_f\right) \ \cong \ H^*_{\text{de Rham}}(X).
$
This immediately yields the strong Morse inequalities.

The restriction of $d$ to $\Stab_f$ has the form
$$
d [S_p] \ = \ \sum_{\lambda_q = \lambda_p-1} n_{pq} [S_q]
$$
where, by a basic result of Federer, the constants
$ n_{pq}$ are integers.  One  concludes that
$
\Stab_f^\BZ  \equiv\ \text{span}_\BZ \left \{ [S_p]\right\}_{p \in {\rm Cr}(f)}
$
is a finite-rank subcomplex of the integral currents  ${\Cal I}_*(X)$
whose inclusion into ${\Cal I}_*(X)$ induces an isomorphism
$
H(\Stab_f^\BZ) \ \cong \ H_* (X;\, \bbz).
$
In the Morse-Smale case the constants $ n_{pq}$
can be computed explicitly by counting flow lines from $p$ to $q$.
This   fact follows directly from Stokes' Theorem  [La].

The operator $\BT$ above can be thought of as the {\it fundamental
solution}  for the de Rham complex provided by the Morse flow.
It is the Morse-theoretic analogue of the Hodge operator $d^*\circ G$
where $G$ is the fundamental solution for the Laplacian (the Greens
operator). The equation (FME) will be referred to as the
 {\it fundamental Morse equation}.  As noted above it naturally
implies the two homology isomorphisms (over $\bbr$ and $\bbz$) which
together encapsule Morse Theory. The proof of this equation employs the
kernel calculus of [HP] to convert current equations on $\xx$ to
operator equations.

The method introduced here has many applications. It was used
in [HL$_2$] to derive a  local version of a formula of MacPherson
[Mac$_{1}$], [Mac$_{2}$] which relates the singularities of a generic bundle map
$A:E\to F$ to characteristic forms of  $E$ and $F$.
However, the techniques apply as well to a significantly larger class of
bundle maps called {\it geometrically atomic}. These include every real
analytic bundle map as well as  generic direct sums and  tensor products
of bundle maps. New explicit formulas relating curvature and
singularities are derived in each case [HL$_3$].  This extends and
simplifies previous work on singular connections and characteristic
currents [HL$_1$].

The approach also works for holomorphic
$\bbc^*$-actions with fixed-points on K\"ahler manifolds.  One finds
a  complex analogue of the current $T$.  General results of
Sommese imply that all the stable and unstable manifolds
of the flow are analytic subvarieties.  One retrieves
 classical results of Bialynicki-Birula [BB] and
Carrell-Lieberman-Sommese [CL], [CS].  The approach also fits
into MacPherson's Grassmann graph construction and  the
construction  of transgression classes  in the refined
Riemann-Roch Theorem [GS].

The method has many other extensions.  It applies to the
multiplication and comultiplication operators in cohomology
whose kernel is the triple diagonal in $\xx\x X$ (\S 11).
This idea can be further elaborated as in [BC].
Similarly, the method fits into constructions of invariants
of knots and 3-manifolds from certain ``Feynman graphs'' (cf.\ [K], [BT]).

The cell decomposition of a manifold by the stable manifolds of a
gradient flow is due to Thom  [T], [S]. That these stable manifolds
embed into the de Rham complex of currents was first observed by
Laudenbach [La].
\medbreak

The authors  are indebted to Janko Latschev for many useful comments.

\demo{Notation}  For an $n$-manifold $X$, $\ce^k(X)$ will
denote the space of smooth $k$-forms on $X$, and ${\cur}^{\lower1pt\hbox{$\scriptstyle k$}}(X)$ the
space of $k$-forms with distribution coefficients (called the {\it currents}
of dimension $n-k$ on $X$).  Under the obvious inclusion
$\ce^k(X)\subset{\cur}^{\lower1pt\hbox{$\scriptstyle k$}}(X)$,  exterior differentiation $d$ has an
extension to all of ${\cur}^{\lower1pt\hbox{$\scriptstyle k$}}(X)$, and $d=(-1)^{k+1}\partial$ where
 $\partial$ denotes the current boundary.
(See [deR] for full definitions.)
We assume here that $X$ is oriented, but this assumption is easily
dropped and all results go through.\enddemo

\section{Finite volume  flows} 

Let $X$ be a compact
oriented manifold of dimension $n$, and  let $\vp_t:X\to X$ be the flow
generated by a smooth vector field $V$ on $X$.  For each $t>0$ consider the
compact submanifold with boundary
$$
\3t\ \equiv \
\{(s,\vp_s(x),x)\, :\, 0\le s\le t\  \text{and}\  x\in X\} \
\subset\ \bbr\times X\times X 
$$
oriented so that 
$$
\bd \3t = \{0\}\x\Delta  - \{t\}\x P_t,
$$
where $\D \subset \xx$ is the diagonal and $P_t =
[\text{graph}(\vp_t)] = \{(\vp_t(x),x)\,:\, x\in X\}$
is the   graph of the diffeomorphism $\vp_t$.
Let $\pr:\bbr\x\xx\to\xx$ denote projection and consider the push-forward
$$
T_t\equiv(\pr)_*(\3t) 
$$
of $\3t$ as a de Rham current.  Since  $\bd$ commutes with
$(\pr)_*$,  we have that
$$
\bd T_t = [\D] - P_t.
$$

Note that the current $T_t$   can be equivalently defined by
$$
T_t \ =\ \Phi_*([0,t]\x X)
$$
where $\Phi : \bbr\x X \to X\x X$ is the smooth mapping defined by
$\Phi(s,x) = (\vp_s(x), x)$.  This mapping is an immersion exactly
on the subset $\bbr\x(X-Z(V))$ where $Z(V) = \{x\in X\,:\, V(x) = 0\}$.
Thus if we fix  a Riemannian metric $g$ on $X$, then $\Phi^*(g\times g)$
is a nonnegative symmetric  tensor which is
$>0$ exactly on the subset  $\bbr\x(X-Z(V))$.

\numbereddemo{Definition} A flow $\1t$ on $X$ is called a {\it finite{\rm -}volume flow} if\break
$\bbr^+\x(X-Z(V))$ has finite volume with respect to the metric
induced by the immersion $\Phi$.
(This concept is independent of the choice of Riemannian
metric on $X$.)
\enddemo

\proclaim{Theorem}  Let $\1t$ be  a finite\/{\rm -}\/volume flow on a compact
manifold $X${\rm ,} and let $T_t$ be the family of currents  defined above with
$$
\bd T_t = [\D] - P_t   \qquad\text{on }\ \xx.
\tag{A}
$$
Then both the limits
$$
P\ \equiv\ \lim_{t\to \infty}P_t \ =\  \lim_{t\to \infty}[\graph\1t]
\qquad\text{and}\qquad
T\ \equiv \  \lim_{t\to \infty}T_t
\tag{B}
$$
exist as currents{\rm ,} and taking the boundary of $T$ gives the
equation of currents
$$
   \partial T  = [\Delta] - P\qquad\text{on }\ \xx.
\tag{C}
$$
\endproclaim

\demo{Proof}
Since $\1t$ is a finite-volume flow, the current $T \equiv
\Phi_*((0,\infty)\x X)$ is the limit in the mass norm of the currents
$T_t = \Phi_*((0,t)\x X)$ as $t \to \infty$.  The
continuity of the boundary operator and equation (A) imply the
existence of  $\lim_{t \to \infty} P_t$ and also establish equation
(C). \enddemo

In the next section we shall see   that Theorem 1.2 has an
important reinterpretation in terms of operators on the differential
forms on $X$.

\numbereddemo{{R}emark}  If we define a relation on $\xx$ by setting $x
\overset{o}\to {\prec} y$ if $y = \1t(x) $ for some $0\leq t<\infty$,
then $T$ is just the (reversed) graph of this relation.  This relation
is always transitive and reflexive, and it is antisymmetric if and
only if $\1t$ has no periodic orbits (i.e., $\overset{o}\to {\prec}$
is a partial ordering precisely when  $\1t$ has no periodic orbits).
\enddemo

\numbereddemo{{R}emark} The immersion $\Phi : \bbr\x (X-Z(V)) \to \xx$
is an embedding outside the subset $\bbr\x {\rm Per}(V)$ where
$$
{\rm Per}(V) = \{x\in X\,:\,  \vp_t(x) = x \ \ \text{for some\ } \ t>0\}
$$
are the nontrivial periodic points of the flow.
Thus, if ${\rm Per}(V)$ has measure zero, then  $T_t$ is given by integration
over the embedded finite-volume submanifold
$\Phi(R_t)$, where $R_t = (0,t)\x \widetilde X$ and
$\widetilde X =  X-Z(V)\cup {\rm Per}(V)$.   If furthermore the flow has finite
volume, then $T$ is given by integration over the embedded, finite-volume
submanifold  $\Phi(R_{\infty})$.

 There is  evidence
that any flow with periodic points cannot have finite volume.
However,  a gradient flow never has periodic points, and many such
flows are  of finite volume.

Note that any flow with fixed points on $S^1$ has finite volume.\enddemo

\numbereddemo{{R}emark}  A standard method for showing that a given flow is
finite volume  can  be outlined as follows.  Pick a
coordinate change $t\mapsto \rho$ which sends $+\infty$ to $0$ and
 $[t_0,\infty]$ to $[0,\rho_0]$.  Then show that
$$
\wh\ct\ \equiv\ \{(\rho, \vp_{t(\rho)}(x),x)\, :\, 0<\rho<\rho_0\}
\ \ \text{has finite volume in} \ \ \bbr\x X\x X.
$$
Pushing forward to $\xx$ then yields the current $T$ with finite
mass.
Perhaps the most natural such coordinate change is $r=1/t$.
Another natural choice (if the flow is considered multiplicatively)
is $s=e^{-t}$.  Of course finite volume in the $r$ coordinate insures
finite volume in the $s$ coordinate since $r \mapsto s=e^{-1/r}$ is a
$C^{\infty}$-map.
\enddemo
 
Many interesting flows can be seen to be finite volume as follows.

\nonumproclaim{Lemma 1.6}  If $X$ is analytic and $\wh\ct \subset \bbr\x\xx$ is
contained in a real analytic subvariety of dimension $n+1${\rm ,} then
$\vp_t$ is a finite\/{\rm -}\/volume flow{\rm .}
\endproclaim

\demo{Proof}
The manifold points of a real analytic subvariety   have
(locally) finite volume. \enddemo

\demo{Example {\rm 1.7}}
  Consider $S^n$ with two coordinate
charts  $\bbr^n$ and the coordinate change $x =y/|y|^2$.  In the
$y$-coordinates consider the translational flow $\vp_t(y) = y+tu$
where $u\in \bbr^n$ is a fixed unit vector. In the $x$-coordinates
this becomes
$$
\vp_t(x) \ =\ { \frac{|x|^2}{|x+t|x|^2u|^2}}
\biggl(x+t|x|^2u\biggr).
$$
To see that $\vp_t$ is finite-volume flow  let
$r=1/t$ and note that  $\wh\ct$
is defined by  algebraic equations
so that Lemma 1.6 is applicable.
The vector field $V$ has only one zero and is therefore
not a gradient vector field.  It is however a limit of gradient vector
fields.
\enddemo

\section{The operator equations}

The current equations (A), (B), (C) of Theorem 1.2 can be translated\break
into operator equations by a kernel calculus which was introduced in [HP]
and which we shall now briefly review.  For clarity of exposition we assume
our manifolds to be orientable.  However, the discussion here extends
straightforwardly to the nonorientable case and  to forms and currents with
values in a flat bundle.

Let $X$ and $Y$ be compact oriented manifolds, and let  $\pi_Y$
and $\pi_X$ denote projection of $\yx$ onto $Y$ and $X$ respectively.
Then each  {current} (or {\it kernel})
$K\in{\cur}^*(Y\x X)$, determines an operator
$\BK:\SE^*(Y)\to{\cur}^*(X)$
by the formula
$$
\BK(\a) = (\pi_X)_* (K\land\piy\a);
$$
i.e., for a  smooth form $\b$ of appropriate degree on $X$
$$
\BK(\a)(\b) = K(\piy\a\land\pix\b).
$$

\demo{Example {\rm 2.1}} Suppose
$\vf:X\to Y$ is a smooth map, and let
$$
P_{\vf} =  \{(\vf(x),x):x\in X\}\subset Y\x X 
$$
be the graph of $\vf$ with orientation induced from $X$.  Then  $\bp\vf$ is
just the
pull-back operator $\bp\vf(\a)\equiv\vf^*(\a)$ on differential forms $\a$.
In particular if $X=Y$ and $K=[\D]\subset\xx$ is the diagonal (the graph of the
identity), then $\BK =\BI:\SE^*(X)\to{\cur}^*(X)$ is the standard inlcusion
of the smooth forms into the currents on $X$.\enddemo

\nonumproclaim{Lemma 2.2}   Suppose an operator $\BT:\SE^*(Y)\to{\cur}^*(X)$
has kernel\break $T\in{\cur}^*(Y\x X)${\rm .} Suppose that $\BT$ lowers degree by
one{\rm ,} or equivalently{\rm ,} that ${\rm deg}(T) =\dim Y-1${\rm .}
Then
$$
\text{The operator }\ d\circ\BT + \BT\circ d\ \text{ has kernel }\
\partial T.\tag2.1
$$
\endproclaim

\demo{Proof}  The boundary operator $\partial$   is the dual of exterior
differentiation. That is $\partial T$ is defined by
$(\partial T)(\piy\a\land\pix \b)=T(d(\piy\a\land\pix \b))$.
Also, by definition, $\BT(d\a)(\b)\equiv
T(\piy(d\a)\land\pix \b)$, and
$$
d(\BT(\a))(\b)\equiv(-1)^{\deg\a} \BT(\a)(d \b) =
(-1)^{\deg\a}T(\piy \a \land\pix d \b),
$$
since $\BT(\a)$ has degree equal to $\deg\a-1$. \enddemo

The results that we need are summarized in the following table.
$$
\matrix
{\rm Operators} &\qquad & {\rm Kernels}\\
\BI& &[\Delta]\\
\pt \equiv\vp^*_t&&  P_t\equiv[\graph\vp_t]\\
\BT&&T\\
d\circ\BT+ \BT\circ d &  & \partial T .
\endmatrix
$$

From this table Theorem  1.2  can be reformulated as follows.

\nonumproclaim{Theorem 2.3}  Let $\1t$ be  a finite\/{\rm -}\/volume flow on a compact
manifold $X${\rm ,} and let $\pt:\SE^*(X)\to\SE^*(X)$   be the operator given
by pull\/{\rm -}\/back
$$
\pt(\a)\equiv\vp^*_{t}(\a).
$$
Let $\bt:\SE^k(X)\to\SE^{k-1}(X)$ be the family of operators associated
to the currents $T_t$ defined in Section~{\rm 1.}
Then for all $t$ one has that
$$
d\circ\bt + \bt\circ d = \BI - \pt. \tag"{(A)}"
$$
Furthermore{\rm ,} the limits
$$
\BT = \lim_{t\to \infty}\bt\quad\text{and}\quad
\BP = \lim_{t\to \infty}\pt\tag"{(B)}"
$$
exist as operators from forms to currents{\rm ,} and they satisfy the
equation
$$
d\circ\BT + \BT\circ d = \BI - \BP.\tag"{(C)}"
$$
\endproclaim

\vglue-12pt
\section{Morse-Stokes gradients; axioms}

Let  $f \in C^{\infty}(X)$ be a Morse function on a compact
$n$-manifold $X$, and let $\crr(f)$ denote the
(finite) set of critical points of $f$.
Recall that   $f$ is  a {\it Morse function} if
its Hessian at each critical point is nondegenerate. The standard
Morse Lemma    asserts that in a neighborhood of each
$p \in \crr(f)$ of index $\l$, there exist {\it
canonical local coordinates} $(u_1,\ldots ,u_{\l}, v_1,\ldots ,v_{n-\l})$
for $|u| < r$, $|v|< r$ with $(u(p),v(p))= (0,0)$ such that
$$
f(u,v) = f(p) - |u|^2 + |v|^2.  \tag{3.1}
$$

  Fix a Riemannian metric $g$ on $X$  and let  $\vf_t$ denote the
flow associated to $\nabla f$.  Suppose that at $p \in \crr(f)$ there
exists a canonical coordinate system $(u,v)$ in which
$g=|du|^2+|dv|^2$. Then  in these coordinates the
 gradient flow is given by
$$
\vf_t(u,v) = (e^{-t} u, e^t v) .  \tag3.2
$$
Metrics $g$ with this property at every $p \in \crr(f)$ will be called
{\it $f$-tame}.

Now to each $p\in\crr(f)$ are
associated the {\it stable} and {\it unstable manifolds} of  the
flow, defined respectively by
$$
S_p = \{x\in X : \lim_{t\to\infty}\vp_t(x) = p\}
\qquad\text{and}\qquad
U_p = \{x\in X : \lim_{t\to -\infty}\vp_t(x) = p\}.
$$
If $g$ is $f$-tame, then for coordinates $(u,v)$ at $p$, chosen
as above, we consider the disks $$
S_p(\e) = \{(u,0) : |u|<\e\}
\qquad\text{and}\qquad
U_p(\e) = \{(0,v) : |v|<\e\}
$$
and observe that
$$
S_p = \bigcup_{-\infty<t< 0}\vp_t\left(S_p(\e)\right)
\qquad\text{and}\qquad
U_p = \bigcup_{0<  t <+\infty}\vp_t\left(U_p(\e)\right).
\tag3.3
$$
Hence, $S_p$ and $U_p$ are contractible submanifolds (but not closed
subsets) of $X$ with
$$
\dim S_p = \la_p \qquad\text{and}\qquad \dim U_p = n-\la_p
\tag3.4
$$
where $\la_p$ is the index of the critical point $p$.
For each $p$ we choose an orientation on $U_p$.   This gives an
orientation on  $S_p$ via the splitting
$T_pX = T_pU_p\oplus T_pS_p$.

The flow $\vp_t$ induces a partial ordering on $X$ by setting
$x \prec y$ if there is a continuous path consisting of a finite
number of foward-time orbits, which begins with $x$ and ends with
$y$. We shall see in the proof of 3.4 that this relation is simply the
closure of the partial ordering of Remark 1.3.

\demo{Definition {\rm 3.1}} The gradient flow of a Morse function $f$ on a
Riemannian  manifold $X$ is called   {\it Morse\/{\rm -}\/Stokes}
if the metric is $f$-tame and if:
\smallbreak
\item{(i)}  The flow is a finite-volume flow.
 \vglue3pt  \item{(ii)}  Each of the stable and unstable manifolds $S_p$ and
$U_p$ for $p\in\crr(f)$ has finite volume.
\vglue3pt    \item{(iii)}  If $p\prec q$ and $p\neq q$ then $\la_p<\la_q$, for all
$p$, $q\in\crr(f)$.
\enddemo

\demo{{R}emark  {\rm 3.2}} If the  gradient flow of
$f$ is Morse-Smale, then it is Morse-Stokes.  Furthermore, for
any Morse function $f$ on a compact manifold $X$ there  exist
riemannnian metrics on $X$ for which the gradient flow of $f$ is
Morse-Stokes. In fact these metric are dense in the space of all
$f$-tame metrics on $X$. See [HL$_4$].
\enddemo

\nonumproclaim{Theorem 3.3}  Let $f \in C^{\infty}(X)$ be a Morse  function on a
compact Riemannian manifold whose gradient flow is Morse\/{\rm -}\/Stokes{\rm .}  Then
the current $P$ in the fundamental equation
$$
\partial T\ =\ [\D] - P    \tag{3.5}
$$
from Theorem 2.1 is given by
 $$
P\ =\ \sum_{p \in {\rm Cr}(f)} [U_p]\times [S_p]    \tag{3.6}
$$
where $U_p$ and $S_p$ are oriented so that $U_p\x S_p$
agrees in orientation with    $X$ at~$p${\rm .}
\endproclaim

\demo{{P}roof} For each critical point $p\in\crr(f)$ we define
$$
\wt U_p\equdef\bigcup_{p\prec q} U_q = \{x\in X:  p\prec  x\}.
$$
\enddemo

\nonumproclaim{Lemma 3.4}   Let $f \in C^{\infty}(X)$ be any Morse function
whose gradient flow is of finite volume{\rm ,} and let  $\spt P\subset\xx$
denote the support of the current $P$ defined in Theorem {\rm 1.2.} Then
$$
  \spt P\subset\bigcup_{p\in\crr(f)}\wt U_p\x S_p.
$$
\endproclaim

\demo{Proof}  Since $P=\lim_{t\to  \infty}P_t$ exists and $P_t=
\{ (\vf_{ t}(x),x)\,:\, x\in X\}$, it is clear
that $(y,x)\in\spt P$ only if there exist sequences $x_i\to x$ in
$X$ and $s_i\to \infty$ in $\bbr$ such that $y_i\equiv \vp_{s_i}(x_i)\to
y$. Let $L(x_i,y_i)$ denote the  oriented flow line from
$x_i$ to $y_i$. Since the lengths  of these lines are
bounded, an elementary compactness argument implies that a subsequence
converges to a piecewise flow line  $L(x,y)$  from $x$ to $y$. By the
continuity of the boundary operator on currents,
$\partial L(x,y) = [y]-[x]$.  Finally, since $s_i\to\infty$ there must be
at least one critical point on $L(x,y)$, and we define
$p=\lim_{s\to\infty}\1s(x)$. \enddemo

Consider now the subset
$$
\Sigma\ \equiv\ \bigcup\lm_{\Sb p\prec q\\ p\neq q\endSb}U_q\x S_p
\ \subset \ \xx 
$$
and define
$
\Sigma' = \Sigma \cup \{(p,p)\,:\,p\in {\rm Cr}(f)\}.
$
From (3.4) and Axiom (iii) in  3.1 note that
$$
\Sigma' \text{ is a finite disjoint union of submanifolds of
dimension} \leq n-1.   \tag3.7
$$

\nonumproclaim{Lemma 3.5}  Consider the embedded submanifold
$$
T \ = \ \{(y,x)\,:\, x \notin  {\rm Cr}(f),\   \text{and $y = \vp_t(x)$
for some } \ 0<t< \infty \} 
$$
of $\xx - \Sigma'${\rm .}  The closure $\overline{T}$ of $T$
in  $\xx - \Sigma'$
is a proper $C^{\infty}$\/{\rm -}\/submanfold with boundary
$$
\bd \overline{T} = \Delta - \sum_{p\in\crr(f)} U_p\x S_p.
$$
\endproclaim

\demo{Proof} We first show that it will suffice  to prove the
assertion in a neighborhood of $(p,p) \in \xx$ for $p\in {\rm Cr}(f)$.
Consider  $(\bar y,\bar x) \in \overline{T} -T$. If  $(\bar y,\bar x)
\notin \Sigma'$, then the proof of Lemma 3.4 shows that either  $(\bar
y,\bar x) \in \D$ or $(\bar y,\bar x) \in  U_p\x S_p$ for some $p\in
{\rm Cr}(f)$.    Near points  $(\bar x,\bar x) \in \D$, $\bar x \notin
{\rm Cr}(f)$,  one  easily checks that  $\overline{T}$ is a submanifold with
boundary $\D$. If $(\bar y,\bar x) \in  U_p\x S_p$, then for
sufficiently large $s>0$, the  diffeomorphism $\psi_s(y,x) \equiv
(\vp_{-s}(y), \vp_s(x))$ will map  $(\bar y,\bar x)$ into any given
neighborhood of $(p,p)$. Note that $\psi_s$ leaves the subset $U_p\x
S_p$ invariant,  and that $\psi_s^{-1}$ maps $T$ into $T$.  Hence, if
$\bd \overline{T} = \D - U_p\x S_p$, in a neighborhood of  $(p,p)$,
then $\bd  \overline{T}  = -U_p\x S_p$ near $(\bar y,\bar x)$.

Now in a neighborhood ${\Cal O}$ of $(p,p)$ we may choose coordinates
as in (3.2) so that $T$ consists of points $(y,x) = (\bar u,\bar
v,u,v)$ with $\bar u =e^{-t}u$ and $\bar v = e^t v$ for some
$0<t<\infty$.  Consequently, in ${\Cal O}$  the set  $\overline T$ is
given by the equations
$$
\bar u =su
\ \ \text{and}\ \
 v = s \bar v
\qquad\text{for some} \ \     0\leq s \leq 1.
$$
This obviously defines a submanifold in ${\Cal O}-\{(p,p)\}$
with boundary consisting of $\D$ and the set $\{\bar u = 0, \ v = 0\}
\cong  (U_p\x S_p)\cap {\Cal O}$. \enddemo

  Lemma 3.5 has the following immediate consequence
$$
\spt\biggl\{P-\sum_{p\in\crr(f)}[U_p]\x [S_p]\biggr\} \subset \Sigma'
.
\tag3.8
$$

We now apply the following elementary but key result of Federer.
\specialnumber{3.6}
\proclaimtitle{[F, 4.1.15]}
\proclaim{Proposition}  Let  $[W]$ be a current in $\bbr^n$
defined by integration over a $k$\/{\rm -}\/dimensional oriented submanifold
$W$ of locally finite volume{\rm .}  Suppose  $\text{\rm supp}\, (d[W])
\subset \bbr^{\ell}${\rm ,} a linear subspace of dimension $\ell < k-1${\rm .}
Then  $d[W]=0${\rm .}
\endproclaim
 
Combining this with (3.7) and (3.8)   proves
(3.6) and completes the proof of Theorem 3.3.  \enddemo

\section{Morse theory}

  Combining Theorems 2.3 and 3.3 gives the following.

\nonumproclaim{Theorem 4.1} Let $f \in C^{\infty}(X)$ be a Morse function on
a compact Riemannian manifold $X$ whose gradient flow $\vf_t$ is
Morse\/{\rm -}\/Stokes{\rm .} Then  for every differential form $\a \in \ce^k(X)${\rm ,}
$0\leq k\leq n${\rm ,} one has
$$  \BP(\a) = \lim_{t\to\infty}\vf_t^*\a  =
\sum_{p\in\crr(f)} r_p(\a) [S_p]
$$
where  the {\rm ``}\/residue\/{\rm ''} $r_p(\a)$ of $\a$ at $p$ is defined by
$
r_p(\a) = \int_{U_p} \a
$
if $k=n-\la$ and $0$ otherwise{\rm .}
Furthermore{\rm ,} there is an operator $\BT$ of degree $-1$ on $\ce^*(X)$ with
values in flat currents{\rm ,} such that
$$
d\circ \BT + \BT \circ d \ =\ \BI - \BP . \tag4.1
$$
\endproclaim

\vglue9pt

Note that
$\BP:\SE^*(X)\lra{\cur}^*(X)$ maps onto the finite-dimensional subspace
$$
\Stab_f^* \ \equdef \ \text{span}_{\bbr}\bigl\{
[S_p]\bigr\}_{p\in\crr(f)},\tag4.2
$$
and that from (4.1) 
$$
\BP\circ d \ =\ d\circ\BP.
$$

\vglue9pt

\nonumproclaim{Theorem 4.2}  The subspace $\Stab_f^*$ is $d$\/{\rm -}\/invariant
and is therefore a subcomplex of ${\cur}^*(X)${\rm .} The surjective
linear map of cochain complexes
$$
\widehat{\BP} : \SE^*(X)\lra\Stab_f^*
$$
defined by $\BP$ induces an isomorphism
$$
\widehat{\BP}_* : H_{\text{deR}}(X) \ \overset{\cong}\to\lra\ H(\Stab^*_f).
$$
\endproclaim

\vglue9pt

\demo{Proof} From Theorem 4.1 it follows immediately that
$${\BP}_*={\BI}_*:H(\SE^*(X))\to H({\cur}^*(X)),$$ and by [deR], ${\BI}_*$ is
an isomorphism. Since  $i_*  \circ \widehat{\BP}_*=\BP_*$, where
$i:\Stab_f^* \to{\cur}^*(X)$ denotes the inclusion, we see that
$\widehat{\BP}_*$ is injective.

However, it does not follow  formally that $\widehat{\BP}_*$ is
surjective  (as  examples show). We must prove that
every $S\in \Stab^*_f$ with $dS=0$ is of the form $S=\BP(\gamma)$
where  $d\gamma=0$.

  To do this we  note
the following consequence of 3.1(iii). For all $p,q\in {\rm Cr}(f)$ of the
same index, one has that $\overline{S}_p\cap \partial
U_q=\emptyset$,  $\overline{S}_p\cap  \overline{U}_q=\emptyset$ if
$p\ne q$, and $\overline{S}_p\cap  \overline{U}_p= \{p\}$.
Now let $S=  \sum r_pS_p$ be a $k$-cycle in $\cs_f$. Then
by the above there is an open neighborhood $N$ of  $\spt \,(S)$  such
that $ N \cap \partial U_q=\emptyset$ for all $q$ of index
$k$.   Since by [deR] cohomology with compact supports can be computed
by either smooth forms or all currents, there exists a
current $\sigma$ with compact support  in $N$  such that $d\sigma =
\gamma-S$ where $\gamma$ is a smooth form. Now for each $q$ of index
$k$, $U_q$ is a closed submanifold of $N$, and we can choose a family
$U_{\e}$ of smooth closed forms in $N$
such that $\lim_{\e\to0}U_{\e}= U_q$  in $\cur(N)$ and
$
\lim_{\e\to 0}(S_p,U_{\e})\ =\ (S_p,U_{q}) \ =\ \delta_{pq}[p]
$
for all $p\in {\rm Cr}(f)$ of index $k$.
This is accomplished via standard Thom form constructions for
the normal bundle using canonical coordinates at  $q$ (e.g.
[HL$_1$]). It follows that
 $(\gamma-S, U_q)= \lim_{\e\to0}(\gamma-S, U_{\e})=
\lim_{\e\to0}(d \sigma, U_{\e})= \lim_{\e\to0}(\sigma, dU_{\e})
=0$. We conclude  that $\int_{U_q}\gamma = (S, U_q)=r_q$ for all
$q$ and so $\BP(\gamma)=S$ as  claimed. \enddemo

The vector space $\Stab_f^*$ has a distinguished lattice
$$
\Stab_f^{\bbz} \ \equdef\ \text{span}_{\bbz}
\left\{[S_p]\right\}_{p \in {\rm Cr}(f)} \ \subset\  \Stab_f
$$
generated by the stable manifolds.  Note that $\Stab_f^{\bbz}$ is a
subgroup of the integral currents  $\SI(X)$ on $X$.

\nonumproclaim{Theorem  4.3}   The lattice $\Stab_f^{\bbz}$ is preserved by
$d${\rm ;} i.e.{\rm ,} $(\Stab_f^{\bbz}, d)$ is a subcomplex of
$(\Stab_f, d)${\rm .}  Furthermore{\rm ,} the inclusion of complexes
$(\Stab_f^{\bbz}, d) \subset (\SI(X), d)$ induces an isomorphism
$$
H(\Stab_f^{\bbz})\ \cong\ H_*(X;\, \bbz)
$$
\endproclaim

\demo{Proof}  Theorem 4.2 implies that for any $p\in {\rm Cr}(f)$ we have
$$
d[S_p] = \sum_{\lambda_q = \lambda_p-1} n_{p,q} \ [S_q]
\tag{4.3}
$$
for real numbers $n_{p,q}$. In particular, $d[S_p]$ has finite mass,
and so by [F; 4.2.16(2)] it is a {\it rectifiable} current.  This
implies that $n_{p,q} \in\bbz$ for all $p,q$, and the first
assertion is proved.

Now the domain of the operator    $\ps$ extends to include any $C^1$
$k$-chain $c$ which is transversal  to the submanifolds $U_p, \  p \in
{\rm Cr}_k(f)$, while the domain of   $\ts$ extends to any
$C^1$ chain $c$ for which $c \times X$   is transversal
to $T$.  Standard
transversality arguments show that such chain groups (over $\bbz$)
compute $H_*(X;\bbz)$. Furthermore they include $\Stab_f^{\bbz}$.
The result then follows from  (4.1).   \enddemo

\nonumproclaim{{C}orollary 4.4}  Let $G$ be a finitely generated abelian group{\rm .}
Then there are natural isomorphisms 
$$
H(\Stab_f^{\bbz}\otimes_{\bbz} G) \ \cong\ H_*(X;G).
$$
\endproclaim

The first assertion of Theorem 4.3 has a completely elementary proof
whenever the flow is  {\it Morse-Smale}, i.e., when $S_p$ is
transversal to $U_q$ for all $p,q \in {\rm Cr}(f)$. Suppose the flow is
Morse-Smale and that $p,q \in {\rm Cr}(f)$ are  critical points with $\l_q
= \l_p -1$.   Then $U_q\cap S_p$ is the  union of a finite set of
flow lines from $q$ to $p$ which we denote $\Gamma_{p,q}$.  To each
$\gamma \in  \Gamma_{p,q}$ we assign an index $n_\g$ as follows.
Let $B_\e \subset S_p$ be a small ball centered at $p$  in a
canonical coordinate system (cf.\ (3.1) ), and  let $y$ be the point
where $\g$ meets $\partial B_\e$.  The orientation of  $S_p$ induces
an orientation on $T_y(\partial B_\e)$, which is identified by
flowing backward along $\g$ with $T_q(S_q)$.  If this identification
preserves orientations we set $n_\g = 1$, and if not, $n_\g = -1$.
We then define
 $$
N_{p,q} \ \equdef \ \sum_{\g \in \Gamma_{p,q}} n_\g.
\tag{4.4}
$$
As in [La] Stokes' Theorem now directly gives us the following (cf.\
[W]).

\nonumproclaim{Proposition 4.5}  When the gradient flow of $f$ is
Morse\/{\rm -}\/Smale{\rm ,} the coefficients in {\rm (4.3)} are given by
$$
n_{p,q} \ =\ (-1)^{\lambda_p}\,N_{p,q}.
$$
\endproclaim

\demo{Proof}  Given a form $\a$ of degree $\lambda_p -1$, we have
$$
(-1)^{\lambda_p}d [S_p](\a) \
= \ \int_{S_p}d\a \ =\ \lim_{r\to \infty} \int_{d S_{p}(r)}  \a
$$
where $S_{p}(r) = \vf_{-r}(S_{p}(\varepsilon))$ as in Section~3.
It suffices to consider forms $\a$ with support near $q$ where
$\la_q = \la_p-1$.  Near such $q$, the set $S_{p}(r)$, for large
$r$,  consists of a finite number of manifolds with boundary,
transversal to $U_q$. There is one for each $\g \in
\Gamma_{p,q}$.   As $r \to \infty$ along one such $\g$, $d
S_{p}(r)$ converges to $\pm S_q$ where the sign is determined by the
agreement (or not) of the orientation of  $d S_{p}(r)$ with the
chosen orientation of    $S_q$.   \enddemo

\demo{{R}emark  {\rm 4.6}} The integers $N_{p,q}$ have a simple definition in
terms of currents.  Set $S_{p}(r) = \vf_{-r}(S_{p}(\varepsilon))$
and  $U_{q}(r) = \vf_{r}(U_{q}(\varepsilon))$  (cf.\ (3.3)). Then
for all $r$ sufficiently large
$$
N_{p,q}\ =\ \int _X[U_{q}(r)]\wedge d\,[S_{p}(r)]
\tag4.5$$
where the integral denotes evaluation on the fundamental class.
\enddemo

Next we observe that our Morse theory enables one to compute
the full integral homology of $X$ (including torsion) using
smooth differential forms.

\demo{Definition {\rm 4.7}} Let $\ce_f^{\bbz}=\BP^{-1}(\cs_f^{\bbz})$  denote the
group of {\it smooth  forms with integral residues}
and  $\ce_f^0\equiv \ker(\BP)$   the subgroup of
{\it  forms with zero residues}.
Note that $\a\in  \ce_f^{\bbz}$ if and only if $\int_{U_p}\a\in \bbz$ for all
$p\in {\rm Cr}(f)$.  Since $\BP$ commutes with $d$, both $\ce_f^{\bbz}$
and $\ce_f^0$ are $d$-invariant, and we have a short exact sequence
of  chain complexes
$$ 0\longrightarrow \ce_f^0\longrightarrow \ce_f^{\bbz}\buildrel{\BP}\over{\longrightarrow} \cs_f^{\bbz}
\longrightarrow 0
\tag4.6
$$
showing that $\ce_f^{\bbz}$ is an extension of the lattice
$\cs_f^{\bbz}$ by the vector space $\ce_f^0$.
\enddemo

\nonumproclaim{Proposition 4.8} The operator $\BP$ induces an isomorphism
$$\BP_*:H(\ce_f^{\bbz}) \buildrel{\cong}\over{\longrightarrow} \ H(\cs_f^{\bbz}).$$
\endproclaim

\demo{Proof} Note that (4.6) embeds into the short exact sequence
$ 0\longrightarrow \ce_f^0 \longrightarrow \ce \buildrel{\BP}\over{\longrightarrow} \cs_f
\longrightarrow 0$, and by 4.2 we know that $\BP_*:H(\ce)\to H(\cs_f)$
is an isomorphism. Hence,
$H(\ce_f^0)=0$ and the result follows. \enddemo

\section{Poincar\'e duality}

In this context there is a simple proof of Poincar\'e duality.
Given two oriented submanifolds $A$ and $B$ of complementary
dimensions in $X$ which meet transversally in a finite number of
points, let $A\bullet B = \int_X[A]\wedge[B]$ denote the algebraic
number of intersections points.  Then for any $k$ we have
$$
U_q\bullet S_p  \ =\ \delta_{pq} \qquad\qquad\text{for all}\ \ p,q
\in {\rm Cr}_k(f)
\tag{5.1}
$$
where ${\rm Cr}_k(f) \equiv \{ p\in {\rm Cr}(f)\ :\ \lambda_p = k\}$. This
gives a formal identification
$$
\Ustab_f^{\bbz}\ \equdef\ \bbz\cdot \bigl\{[U_p]\bigr\}_{p\in {\rm Cr}(f)}
\ \cong \  \text{Hom}\left(\Stab_f^{\bbz},\, \bbz\right).
\tag{5.2}
$$
Therefore, taking the adjoint of $d$ gives a differential $\d$
on $\Ustab_f^{\bbz}$ with the property that
$H_{n {-} *}(\Ustab_f^{\bbz}, \delta) \cong H^*(X; \bbz)$.
On the other hand the arguments of Sections 1--4    (with $f$
replaced by $-f$) show that  $ \Ustab_f^{\bbz}$ is $d$-invariant
with  $H_*(\Ustab_f^{\bbz}, d) \cong H_*(X; \bbz)$.   However, these
two differentials  on  $\Ustab_f^{\bbz}$ agree up to sign as we see
in the next lemma.

\nonumproclaim{Lemma 5.1}  One has
$$
(dU_q)\bullet S_p\ =\ (-1)^{n-k} \, U_q\bullet (dS_p)
\tag{5.3}
$$
for all $p \in {\rm Cr}_k(f)$ and $q \in {\rm Cr}_{k-1}(f)${\rm ,} and for any $k${\rm .}
\endproclaim

\demo{Proof}
One can see directly  from the definition that the integers $N_{p,q}$
are invariant (up to a global sign) under time-reversal in the flow.
However, for a simple current-theoretic proof consider the
1-dimensional current $[U_q(r)]\wedge [S_p(r)]$ consisting of a
finite sum of oriented line-segments in the flow lines of
$\Gamma_{p,q}$ (cf.\ Remark  4.6). Note that
$$
d\, [U_q(r)]\wedge [S_p(r)]\ =\
(d\,[U_q(r)])\wedge [S_p(r)] \ +\
(-1)^{n-k+1}[U_q(r)]\wedge (d\,[S_p(r)])
$$
and apply (4.5).\enddemo

\nonumproclaim{{C}orollary 5.2 {\rm (Poincar\'e Duality)}}
$$
H^{n-k}(X;\bbz) \ \cong\ H_k(X;\bbz) \qquad\text{{\it for all}} \ \
k. $$
\endproclaim

\demo{Note {\rm 5.3}} The Poincar\'e duality isomorphism can be
realized in our operator picture as follows. Consider the ``total
graph'' of the flow: in $\xx$
$$
T_{\text{tot}}\ =\ T^*+T \ =\
\{(y,x)\,:\, y=\vp_t(x)\ \ \text{for some}\ t\in \bbr\}
$$
with corresponding operator ${\BT}_{\text{tot}}$.  Then one has
the operator equation
$$
d\circ  {\BT}_{\text{tot}} +  {\BT}_{\text{tot}}\circ d\ =
\  \BP - \check{\BP},
$$
where
$$
\BP  = \sum_{p\in {\rm Cr}(f)} [U_p]\times[S_p] \qquad\text{and} \qquad
\check{\BP}  = \sum_{p\in {\rm Cr}(f)} [S_p]\times[U_p].
$$
This chain homotopy induces an isomorphism
$H_*(\Ustab_f^{\bbz}) \cong H_*( {\Stab}_f^{\bbz})$, which after
identifying $\Ustab_f^{\bbz}$  with the cochain complex via (5.1)
and (5.3), gives the duality isomorphism 5.2.
When $X$ is not oriented, a parallel analysis yields Poincar\'e
duality with mod 2 coefficients.
\enddemo

\section{Generalizations}  As seen in Section~2, for
any  flow  $\vf_t$ of finite volume the operator $\ps(\a) = \lim_{t\to
\infty}\vf_t^*(\a)$ exists and is chain homotopic to the identity.
This situation occurs often.  Suppose for example that
 $f:X\to \bbr$ is a smooth function whose critical set is a finite
disjoint union $$
{\rm Cr}(f) \ = \ \coprod_{j=1}^\n F_j
$$
of compact submanifolds $F_j$ in $X$ and   that Hess$(f)$ is  
nondegenerate  on the normal  spaces to ${\rm Cr}(f)$.  Then for any
$f$-tame gradient flow $\vf_t$ (cf.\ [L]) there are stable and unstable
manifolds $$
S_j = \{x\in X \ :\ \lim_{t\to\infty} \vf_t(x)  \in F_j\}
\qquad\text{and}\quad
U_j = \{x\in X \ :\ \lim_{t\to -\infty} \vf_t(x)  \in F_j\}
$$
for each $j$, with projections
$$
S_j \buildrel{\tau_j}\over\longrightarrow F_j \buildrel{\sigma_j}\over\longleftarrow U_j,
\tag{6.1}
$$
where $$
\tau_j(x) = \lim_{t\to\infty} \vf_t(x) \qquad\text{and}\qquad
 \sigma_j(x) = \lim_{t\to -\infty} \vf_t(x).
$$
For each $j$, let $n_j = \text{dim}(F_j)$ and set $\lambda_j =
\text{dim}(S_j) - n_j$.  Then dim$(U_j) = n-\lambda_j$. For $p \in
F_j$ we define $\lambda_p \equiv\lambda_j$ and $n_p \equiv n_j$.

\demo{Definition {\rm 6.1}} The gradient flow $\vf_t$ of a smooth function $f \in \Ci(X)$
on a Riemannian manifold $X$ is called a  {\it generalized Morse-Stokes
flow} if is $f$-tame and
\medbreak
\item{(i)} The critical set of $f$ consists of a finite number of
submanifolds $F_1,\ldots ,F_\n$ on the normals of which Hess$(f)$ is
nondegenerate.

\smallbreak\item{(ii)} The manifolds $T$, and $T^*$, and   the stable and unstable
manifolds $S_j$, $U_j$ for $1\leq j\leq \n$ are submanifolds of finite
volume. Furthermore, for each $j$, the fibres of the projections
$\tau_j$ and $\sigma_j$ are of uniformly bounded volume.

\smallbreak \item{(iii)} $p \prec q \ \ \ \Rightarrow
\ \ \ \l_p + n_p < \l_q$ for all $p,q \in {\rm Cr}(f)$.
\enddemo

\nonumproclaim{Theorem  6.2} Suppose $\vf_t$ is a gradient flow  satisfying
the  generalized Morse\/{\rm -}\/Stokes conditions {\rm 6.1} on a compact oriented
manifold $X${\rm .}  Then there is an equation of currents
$$
\partial T \ = \ [\D] - P
$$
 on $X\times X${\rm ,} where $T${\rm ,} $\D${\rm ,} and $P$ are as in
Theorem {\rm 1.2,} and
$$
P \ = \ \sum_{j=1}^\n \left[U_j\times_{F_j}S_j\right]
\tag {6.2}
$$
where $U_j\times_{F_j}S_j \equiv \{(y,x) \in U_j\times S_j \subset
X\times X\ :\ \sigma_j(y) = \tau_j(x) \}$ denotes the fibre product of
the projections {\rm (6.1).}
\endproclaim

\demo{Proof}  The argument  follows closely the proof of Theorem 3.3.
Details are omitted.   \enddemo

In [L] Janko Latschev has found a Smale-type condition which yields this
result in cases where when 6.1 (iii) fails.
In particular, his condition implies only that: $p\prec q
\Rightarrow \lambda_p < \lambda_q$.

 The result above can be translated into the following operator form.

\nonumproclaim{Theorem  6.3} Let $\vf_t$ be a gradient flow  satisfying
the  generalized Morse\/{\rm -}\/Stokes conditions on a
manifold $X$ as above{\rm .}    Then for all smooth forms
$\a$ on $X${\rm ,} the limit
$$
\ps(\a) = \lim_{t \to \infty}\vf_t^*(\a)
$$
exists and defines a continuous linear operator
$\ps : \ce^*(X) \arr {\cur}^*(X)$ with values in
flat currents on
$X${\rm .} This operator fits into a chain homotopy
$$
d\circ \ts + \ts\circ d\ =\BI - \ps.
\tag{6.3}
$$
Furthermore{\rm ,} $\ps$ is given by the formula
$$
\ps(\a) \ =\ \sum_{j=1}^\n {\rm Res}_j(\a) \left[S_j\right]
\tag{6.4}
$$
where
$$
{\rm Res}_j(\a)\ \equiv\
\tau_j^*\left\{(\sigma_j)_*\left(\a\bigl|_{U_j}\right)\right\}.
 \tag{6.5}
$$
\endproclaim

\vglue9pt

\demo{Proof}  This is a direct consequence of Theorem 6.2  except for the
formulae (6.4)--(6.5).  To see this consider the pull-back square
$$\matrix \noalign{\vskip4pt}
U_j\times_{F_j}S_j &\buildrel{t_j}\over\lrar & U_j \\
\noalign{\vskip4pt}
 {\scriptstyle s_j}\Big\downarrow\quad&&     \quad \Big\downarrow {\scriptstyle\sigma_j}  \\
\noalign{\vskip4pt}
S_j &\buildrel{\tau_j}\over\lrar&     F_j\\
\noalign{\vskip4pt}
\endmatrix\tag{6.6}
$$
where $t_j$ and $s_j$ are the obvious projections.  From the
definitions  one sees that
$$
\ps(\a) \ = \ \sum_{j=1}^\n
(s_j)_*\left\{ (t_j)^* \left( \a\bigl|_{U_j} \right) \right\}.
 $$
The commutativity of the diagram (6.6) allows us to rewrite these
terms as in (6.5). \enddemo

\nonumproclaim{{C}orollary 6.4}   Suppose that $\lambda_p + n_p + 1 < \lambda_q$ for
all critical points $p\prec q${\rm .} Then the homology of $X$ is spanned by
the images of the groups $H_{\lambda_j + \ell}(\overline{S_j} )$  for
$j = 1,\ldots ,\n$ and $\ell \geq 0${\rm .}
\endproclaim

\demo{Proof} Under this hypothesis $\partial (U_j\times_{F_j}S_j) = 0$ for all
$j$, and so (6.2) yields a decomposition of $\BP$ into operators
that commute with $d$. \enddemo

  Each map $\tau_j: S_j \to F_j$  can be given the structure of a
vector bundle of rank $\lambda_j$.    The closure $\overline{S_j}
\subset X$ is a compactification of this bundle with a complicated
structure at infinity  (cf.\ [CJS]).  \ There is nevertheless a
homomorphism  $\Theta_j :H_*(F_j) \arr
H_{\lambda_j+*}(\overline{S_j})$ which after pushing forward to the
one-point compactification of ${S_j}$, is the Thom isomorphism.
This leads to the following (cf.\ [AB]).

\nonumproclaim{{C}orollary 6.5} Suppose that $\lambda_p + n_p + 1 < \lambda_q$ for
all critical points $p\prec q$ and that $X$ and all $F_j$ and $S_j$ are
oriented{\rm .}  Then there is an isomorphism 
$$
H_*(X) \ \cong\ \bigoplus_j H_{*-\lambda_j}(F_j)
$$
\endproclaim

One can drop the orientation assumptions by taking
homology with appropriately twisted coefficients.  Extensions of
Theorem 6.3 and Corollary 6.5 to integral homology groups    are
found in [L].      Latschev also derives a spectral sequence, with
geometrically computable differentials, associated to any {\it
Bott-Smale} function satisfying a certain transversality hypothesis.
Assuming  that everything is oriented, the $E^1$-term  is given by $
E^1_{p,q} \ = \ \bigoplus_{\lambda_j= p} H_{q}(F_j; \bbz)  $
and $E^k_{p,q} \Rightarrow H_*(X; \bbz).$

\demo{{R}emark  {\rm 6.6}}  In standard Morse Theory one  considers a proper
exhaustion function $f:X\to \bbr$ and  studies the change in
 topology as one passes from $\{x : f(x) \leq a\}$ to $\{x : f(x)
\leq b\}$.  Our operator approach is easily adapted to this case
(See [HL$_4$]).
\enddemo

\section{The local MacPherson formula and generalizations}  

The ideas
above can be applied to the study of curvature and singularities. Suppose
$$
\a : E \arr F
$$
is a map between smooth vector bundles with connection over a
manifold $X$.  Let $G = G_k(E\oplus F) \to X$ denote the Grassmann
bundle of  $k$-planes in $E\oplus F$ where $k = \text{rank}(E)$.
There is a flow $\vf_t$ on $G$ induced by the flow
$\psi_t :E\oplus F \to E\oplus F$ where $\psi_t(e,f) = (te,f)$.
 On the ``affine chart" Hom$(E,F)$, this flow has the form $\vf_t(A) =
\tfrac 1 t A$.
(Note that here it is more natural to consider the flow
multiplicatively ($\vp_{ts} = \vp_t \circ \vp_s$) than additively.)

\demo{Definition {\rm 7.1}} The section $\a$ is said to be {\it geometrically atomic}
if its {\it radial span}
$$
\ct(\a) \equdef \left\{\tfrac 1 t \a(x) \, : \, 0<t \leq 1
\ \ \text{and}\ \  x\in X\right\}
$$
has finite volume in $G$.
\enddemo

This hypothesis is sufficient to guarantee the existence of
$\lim_{t\to 0} \a_t^*\Phi$ where $\a_t \equiv \frac 1 t \a$ and where
$\Phi$ is   any differential form on
$G$.  Choosing $\Phi = \Phi_0(\Omega_U)$ where  $\Phi_0$ is an
Ad-invariant polynomial on  ${\hbox{\fr g}\hbox{\fr l}}_k(\bbr)$ and
$\Omega_U$ is the curvature of the tautological $k$-plane bundle over
$G$, one obtains formulas on $X$ which relate the singularities of $\a$
to characteristic forms of $E$ and $F$. For {\it normal maps}, maps
which are transversal to the universal singularity sets, this
process yields a local version of  a classical formula of MacPherson
[Mac$_{1}$], [Mac$_{2}$]. In other cases, such as generic direct sum or
tensor products mappings, the process yields new formulas.
The method also applies whenever $\a$ is real analytic.  Details
appear in [HL$_{2}$], [HL$_{3}$].

\section{Equivariant Morse theory} 

 Our approach carries
over to equivariant cohomology by using
Cartan's equivariant de Rham theory (cf.\ [BGV]).
Suppose $G$ is a compact Lie group with Lie algebra $\gg$ acting
on a compact $n$-manifold $X$.  By an {\it equivariant differential
form} on $X$ we mean a $G$-equivariant polynomial map $\a:\gg \arr
\ce^*(X)$. The set of equivariant forms is denoted by
$$
\ce^*_G(X)\ =\ \{S^*(\gg^*)\otimes \ce^*(X)\}^G
$$
and is graded  by declaring elements of $S^p(\gg^*)\otimes \ce^q(X)$
to have total degree $2p+q$.  The equivariant differential
$d_G : \ce^*_G(X) \arr \ce^{*+1}_G(X)$ is defined by setting
$$
\left(d_G\a\right)(V)  \ = \ d\a (V) - i_{\tv}\a(V)
$$
for $V\in \gg$ where  $i_{\tv}$ denotes contraction with the vector field
$\tv$ on $X$ corresponding to $V$.
The complex ${\Cal D}^{'*}_G(X)$ of {\it equivariant currents} is similarly
defined by replacing smooth forms $\ce^*(X)$ by forms ${\Cal D}^{'*}(X)$ with
distribution coefficients.

Consider  a $G$-invariant function $f \in C^\infty(X)$ and suppose
 $X$ is given   a $G$-invariant Riemannian metric.
Suppose  the ($G$-invariant) gradient flow has finite volume and let
$$
\partial T \ = \ \D - P
\tag{8.1}
$$
be the current equation derived in Section~3.  Let $G$ act on $\xx$ by
 $g\cdot (x,y) = (gx, gy)$.

\phantom{odd}

\nonumproclaim{Lemma 8.1}  The current $T$ satisfies\/{\rm :}
\medbreak
\item{\rm (i)} $g_*T = T$ \qquad for all $g\in G${\rm ,} and
\medbreak\item{\rm  (ii)} $i_{\tv}T = 0$ \qquad for all $V\in \gg${\rm .}
\medbreak
\noindent
So also does the current $P${\rm .}
\endproclaim

\vglue6pt

\demo{Proof}
The current $T$ corresponds to integration over the finite-volume
submanifold $\{(x,y) \in \xx-\D\ :\  y = \vf_t(x) \ \ \text{for some\ }
0<t<\infty\}$.
Since $g\vf_t(x) = \vf_t(gx)$, assertion (i) is clear.  The
invariance of $T$ implies the invariance of $P$ by (8.1). For
assertion (ii) note that for any $(n+2)$-form $\beta$ on $\xx$
$$
\left(i_{\tv}T \right) (\beta) \ = \ \int_T i_{\tv}\beta \ =\ 0
$$
since $V$ is tangent to $T$.  Similarly, $i_{\tv}\D = 0$.
Since $d \circ i_{\tv} + i_{\tv}\circ d = {\Cal L}_{\tv}$ (Lie
derivative) and ${\Cal L}_{\tv} T = 0$, we conclude that
$i_{\tv} P = 0$.  \enddemo

\nonumproclaim{{C}orollary 8.2} Consider $T \equiv 1 \otimes T \in 1 \otimes
{\Cal D}^{'n-1}(\xx)^G$ as an  equivariant current of total degree $(n-1)$
on $\xx${\rm .}  Consider $\D$ and $P$ similarly as equivariant currents of
 degree $n${\rm .}  Then
$$
\partial_G T \ = \D - P  \qquad \text{on\ \ }\xx.
\tag{8.2}
$$
\endproclaim

The correspondence between operators and
kernels discussed in\break Section~2 carries over to the equivariant
context.  Currents in ${\Cal D}^{'n-\ell}_G(\xx)$ yield 
$G$-equivariant operators $\ce^*(X) \arr {\Cal D}^{'*+\ell}(X)$; hence
 operators $\ce^*_G(X)  \arr {\Cal D}^{'*+\ell}_G(X)$. Equations of type
(8.2) translate into operator equations
$$
d_G \circ \ts + \ts \circ d_G \ = I - \ps.
\tag{8.3}
$$
Applying the arguments of Section~4 proves the following.

\nonumproclaim{Proposition 8.3}    Let $\vf_t$ be an invariant flow on a compact
$G$\/{\rm -}\/manifold $X${\rm .}  If $\vf_t$ has finite volume{\rm ,} then the limit
$$
\ps (\a) \ = \ \lim_{t\to \infty}  \vf_t^*\a
\tag8.4
$$
exists for all $\a \in \ce^*_G(X)$ and defines a continuous linear
operator $\ps : \ce^*_G(X)\! \arr {\cur}^*_G(X)$ of degree $0${\rm ,} which is
equivariantly chain homotopic to the identity on $\ce^*_G(X)${\rm .}
\endproclaim

The  fact that  $S^*(\gg^*)^G
\cong H^*(BG)$  leads as in Section~4 to  the following result.

\nonumproclaim{Theorem  8.4} Let $f \in  C^\infty(X)$ be an invariant Morse
function on a compact riemnannian $G$\/{\rm -}\/manifold whose gradient flow
$\vf_t$ is Morse\/{\rm -}\/Stokes{\rm .}  Then the continuous linear operator {\rm (8.4)}
defines a map of equivariant complexes
$$
\ps : \ce^*_G(X) \arr S^*(\gg^*)^G\otimes \Stab_f
$$
where $\Stab_f = \text{\rm span}\{[S_p]\}_{p\in {\rm Cr}(f)}$ as in {\rm (4.2)}
and where the differential on\break $S^*(\gg^*)^G\otimes \Stab_f$ is
$1\otimes \partial${\rm .}  This map induces an isomorphism $\phantom{\sum^1}$
$$
H^*_G(X) \bra{\cong}  H^*(BG)\otimes H^*(X).
$$
\endproclaim

Examples of this phenomenon arise in moment map constructions.
For a simple example consider $G = (S^1)^{n+1}/\D$ acting on
$\bbp^n_\bbc$ via the standard action on homogeneous coordinates
$[z_0,\ldots ,z_n]$, and set $f([z]) = \sum k|z_k|^2/\|z\|^2$.
One sees immediately the well-known fact that $H^*_G(\bbp^n_\bbc)$
is a free $H^*(BG)$-module with one generator in each dimension
$2k$ for $k=0,\ldots ,n$.  This extends to all generalized flag manifolds
and to products.

 J. Latschev has pointed out that  there exists an
invariant Morse function for which no choice  of invariant metric
gives a Morse-Stokes flow.

However, the method applies to quite general functions
and yields results as in Sections~6 and 7.  Suppose for example that $f$ is an
invariant function whose critical set consists of a finite number of
nondegenerate critical orbits
$\co_i= G/H_i$, $i=1,\ldots ,N$\  (the generic case).
There is a spectral sequence
with (assuming for simplicity that everything is oriented)
$$
E^1_{p,*}\ = \ \bigoplus_{\lambda_i = p} H^*_G(\co_i)
\ = \ \bigoplus_{\lambda_i = p} H^*(BH_i)
$$
and computable differentials such that
$
E^k_{*,*}\ \Rightarrow\ H^*_G(X)
$
 (cf.\ [L]).

\section{Holomorphic flows and the Carrell-Lieberman-Sommese
theorem}

This method also applies to the
holomorphic case.  For example, given a  $\bbc^*$-action $\vf_t$ on a
compact \Kahler manifold $X$, there is a complex graph
$$
\cct \ \equdef\ \{(t, \vf_t(x),x) \in \bbc^*\times X\times X \, :\,
t \in   \bbc^*\ \ \text{and}\ \ x \in X\}
\ \subset \ \bbp^1(\bbc) \times X\times X
$$
analogous to the graphs considered above.
The following is a result of Sommese [So].

\nonumproclaim{Theorem  9.1} If $\vf_t$ has fixed\/{\rm -}\/points{\rm ,} then
$\cct$ has finite volume and its closure $\overline \cct$
in $\bbp^1(\bbc) \times X\times X$ is an analytic subvariety{\rm .}
\endproclaim

The relation of $\bbc^*$-actions to Morse-Theory is classical.
The action $\vf_t$ decomposes into an ``angular'' $S^1$-action and a
radial flow.  Averaging a \Kahler metric over   $S^1$  and applying
an argument of Frankel [Fr], we find a function $f:X \to \bbr$
of Bott-Morse type  whose gradient  generates the radial action.
Theorem 9.1 implies that the associated current $T$ on $\xx$, which is
a real slice of $\cct$, is real analytic and hence of finite volume.
One can then apply the methods of this paper, in particular those of
Section~6.

 Thus when $\vf_t$ has fixed-points, $\overline \cct$ gives a
rational equivalence between the diagonal $\D$ in $X\times X$ and an
analytic cycle $P$ whose components consist of fibre products of
stable and unstable manifolds over components of the fixed-point set
of the action.

When the fixed-points are all isolated, $P$ becomes a
sum of analytic K\"unneth components
$P =  \sum {\ov S}_p\times {\ov {\widetilde U}}_p$, and we recover the
well-known fact that the  cohomology of $X$ is freely generated by the
stable subvarieties $\{{\ov S}_p\!\}_{p \in \text{Zero}(\vf)}$.  It
follows that $X$ is algebraic and that all cohomology theories
on $X$ (eg. algebraic cycles modulo rational equivalence, algebraic
cycles modulo algebraic equivalence, singular cohomology) are naturally
isomorphic. (See [BB], [ES], [Fr].)

When the fixed-point set  has positive dimension, one
can recover results of Carrell-Lieberman-Sommese for $\bbc^*$-actions
([CL], [CS]), which assert among other things that if
$\dim(X^{\bbc^*}) = k$, then $H^{p,q}(X) = 0$ for $|p-q| > k$.

\section{Local coefficient systems}

    Our method applies
immediately to forms with coefficients in a flat bundle $E\to X$.
In this case the kernels of Section~2 are currents on $X\x X$ with coefficients
in $\Hom(\pi_1^*E,\pi_2^*E)$.  Given a Morse-Stokes flow $\vf_t$ on $X$
we consider the kernel $T_E = h\otimes T$ where $T$ is defined as in Section~3
and $h:E_{\vf_t(x)} \to E_x$ is parallel translation along the flow line.
One obtains the equation $\partial T_E = \Delta_E - P_E$ where $\D_E =
\text{Id}\otimes \D$ and $P_E = \sum_p h_p\otimes ([U_p]\x [S_p])$ with
$h_p:E_y\to E_x$ given by parallel translation along the broken flow line.
Thus, $h_p$ corresponds to \ $\text{Id}:E_p\to E_p$ under the canonical
trivializations $E\bigl|_{U_p} \cong U_p\x E_p$ and
$E\bigl|_{S_p} \cong S_p\x E_p$. We obtain the operator equation
$$
d\circ \BT_E + \BT_E \circ d \ =\ \BI - \BP_E ,\tag10.1
$$
where $\BP_E$ maps onto the finite complex
$$
\SS_E\ \equdef\ \bigoplus_{p\in {\rm Cr}(f)} E_p  \otimes [S_p] 
$$
by integration of forms over the unstable manifolds.  The restriction of
$d$ to $\SS_E$ is given as in (4.3) by
$
d(e\otimes [S_p])  = \sum h_{p,q}(e) [S_q]
$
where $h_{p,q} = (-1)^{\lambda_p}\sum_{\gamma}  h_{\gamma}$ and
$h_{\gamma}:E_p\to E_q$ is parallel translation along $\gamma \in
\Gamma_{p,q}$.  By (10.1) the complex $(\SS_E,d)$ computes $H^*(X;\, E)$.

Reversing time in the flow shows that the complex $\Ustab_E =
\bigoplus_{p} E_p^*  \otimes [U_p]$ with differential defined as
above computes $H^*(X;\, E^*)$. As in Section~5 the dual pairing of
these complexes establishes  Poincar\'e duality.
Furthermore, this  extends to integral currents twisted
by representations of  $\pi_1(X)$ in $\text{GL}_n(\BZ)$ or
$\text{GL}_n(\BZ/p\BZ)$ and gives duality with local coefficient
systems.

\def\bwedge{{\wedge}}

\section{Cohomology operations}

   Our method has many
 extensions.  For example, consider the triple
diagonal $\D_3 \subset X\x X\x X$ as the kernel of the wedge-product
operator \linebreak
$
\bwedge:\SE^*(X)\otimes \SE^*(X)\to\SE^*(X).
$
Let $f$ and $f'$ be functions with Morse-Stokes flows $\vf_t$ and
$\vf'_t$ respectively. Assume that for all
$(p,p')\in {\rm Cr}(f)\x {\rm Cr}(f')$  the stable manifolds $S_p$ and
$S'_{p'}$ intersect transversely in a manifold of finite
volume, and similarly for the unstable manifolds
$U_p$ and $U'_{p'}$.  Degenerating
$\D_3$ gives a kernel $$
T\ \equiv\ \{(\vf_t(x), \vf'_t(x), x) \in X\x X\x X\,:\,
x\in X \ \ \text{and}\ \ 0\leq t <\infty\}
$$
and a corresponding operator $\BT:\SE^*(X)\otimes
\SE^*(X)\to{\Cal D}^{'*}(X)$ of degree -1.  One calculates that
$
\partial T \ =\ \D_3 - M
$
where
$
M = \sum_{(p,p')} \, [U_p]\x [U'_{p'}]\x
[S_p\cap S'_{p'}].
$
The corresponding operator $\BM:\SE^*(X)\otimes \SE^*(X)\to{\Cal D}^{'*}(X)$
is given by
$$
\BM(\a,\b)\ =\ \sum_{(p,p')\in {\rm Cr}(f)\x {\rm Cr}(f')}
\left(\int_{U_p}\a\right)\left(\int_{U'_{p'}}\b\right)
[S_p\cap S'_{p'}].
\tag11.1
$$
The arguments of Sections~1--3 adapt to prove the following.

\nonumproclaim{Theorem  11.1}  There is an equation of operators
$\bwedge - \BM = d\circ\BT+\BT\circ d$ from $\ce^*(X\x X)$ to
${\Cal D}^{'*}(X)$ {\rm (}\/where $\bwedge$ denotes restriction to the diagonal\/{\rm ).}
In particular for $\a, \b \in \ce^*(X)$ we have the chain homotopy
$$
\a\bwedge\b - \BM(\a,\b)\ =\ d\BT(\a,\b) +\BT(d\a,\b) +
(-1)^{\text{deg}\a} \BT(\a,d\b).
\tag11.2
$$
\endproclaim

Note that the operator $\BM$ has range in the finite-dimensional
vector space
$
\SM  \equdef$  $\text{span}_{\bbr}
\{\ [S_p\cap S'_{p'}]\}_{_{(p,p')}}.
$
It converts a pair of smooth forms $\a, \b$ into a linear
combination of the pairwise intersections of the stable manifolds
$[S_p]$ and $[S'_{p'}]$.  If $d\a = d\b = 0$, then
$
\BM(\a,\b)\ =\ \a\wedge\b-d\BT(\a,\b),
$
and so $\BM(\a,\b)$ is a cycle homologous to the wedge product
$\a\wedge\b$. This  operator maps onto the
subspace  $\SM$, and for forms
$\a, \b \in \ce^*(X)$, it satisfies the equation
$$
d\BM(\a, \b)\ =\ \BM(d\a,\b) + (-1)^{\text{deg}\a} \BM(\a,d\b).
\tag11.3
$$

It follows that $d(\SM)\subset \SM$.  Furthermore, for
$(p,p')\in {\rm Cr}(f)\x {\rm Cr}(f')$ one has
$$
d\,[S_p\cap S'_{p'}]\ =\ \sum_{ q\in {\rm Cr}(f)}
n_{pq} [S_q\cap S'_{p'}] +(-1)^{n-\l_p}
\sum_{ q'\in {\rm Cr}(f')}
n'_{p'q'} [S_p\cap S'_{q'}]
\tag11.4
$$
where the $n_{pq}$ are defined as in Section~4.  Thus we retrieve the
cup product over the integers in the Morse complex.

A basic example of a pair satisfying these
hypotheses is simply $f, -f$ where the gradient flow is
Morse-Smale.   Here $U'_p = S_p$
and $S'_p = U_p$ for all $p\in {\rm Cr}(f)$, and we have:

\nonumproclaim{Proposition 11.2} Suppose the  gradient flow of $f$  is
Morse\/{\rm -}\/Smale and $Y$  is a cycle in $X$ which is transversal
to    $U_p\cap S_{p'}$ for all $p,p'\in {\rm Cr}(f)${\rm .} Then for any closed
forms $\a, \b$
 $$
\int_Y \a\wedge\b \ = \ \sum_{p,p'\in {\rm Cr}(f)}
\left(\int_{U_p}\a\right)\left(\int_{S_{p'}}\b\right)
[S_p\cap U_{p'}\cap Y].
$$
\endproclaim

Analogous formulas are derived by the same method for the
comultiplication operator.
Following ideas in [BC] one can extend these constructions to other
cohomology operations.

\section{Knot invariants}

  From the finite-volume flow 1.7 our
methods construct an $(n+1)$-current $T$ on $S^n\x S^n$ with the
property that $
\partial T \ =\  S^n\x \{*\} + \{*\}\x S^n.
$
This is a singular analogue of the form used by Bott and Taubes [BT],
and can be used to develop a combinatorial version of their study of
knot invariants.

\section{Differential characters}

  The ring of differential
characters on a smooth manifold, introduced by Cheeger and Simons [ChS]
in 1973, has played an important role in geometry. In de~Rham-Federer
formulations of the theory (see [HL$_5$]), differential characters are
represented by {\it sparks}.  These are currents $T$ with the property
that $dT =\a-R$ where $\a$ is  smooth and $R$ is  integrally flat.
Our Morse Theory systematically produces such creatures.

Consider the group $\cz^{\bbz}_f= \ce^{\bbz}_f\cap \ker d$  of
closed forms with integral residues; (see \S 4). Then
$\BT:\cz^{\bbz}_f \to \cur$ is a continuous map with the property
that   $
d\BT(\a) \ =\ \a-\BP(\a)
$
for all $\a$. Thus $\BT$ gives a continuous map of $\cz^{\bbz}_f$
into the space of sparks, and therefore also into the ring of
differential characters on $X$.

\input Harvey.ref
\bye

\item {[AB]} M. Atiyah and R. Bott, {\it The Yang-Mills equations
over a Riemann surface},  Phil. Trans., R. Soc. London. A   {\bf 308}
(1982), 523-615.

\item {[BC] }  M. Betz and R.L. Cohen, {\sl Graph  moduli spaces and
cohomology operations }, Stanford Preprint, 1993.

\item {[BB]} A. Bialynicki-Birula, {\sl Some theorems on actions of
algebraic groups},  Annals of  Math., {\bf 98} (1973), 480-497.

\item {[BGV] }  N. Berline, E. Getzler and M. Vergne, {\sl
Heat Kernels and the Dirac Operator}, Grundl. der Math. Wiss. Band 298,
Springer, Berlin-Heidelberg-New York, 1992.

\item {[BT] }  R. Bott and C. Taubes, {\sl On the self-linking of
knots}, J. Math. Phys., {\bf 35 }(10) (1994), 5247-5287.

\item {[CL]} J. B.  Carrell and D. I.  Lieberman, {\sl Holomorphic vector
manifolds and compact \Kahler manifolds},  Invent.  Math., {\bf 21} (1973),
303-309.

\item {[CS]} J. B.  Carrell and A. J. Sommese, {\sl Some topological
aspects of ${\scriptstyle\bbc}^*$-actions on compact \Kahler manifolds},  Comm.  Math.
Helv., {\bf 54} (1979), 583-800.

\item {[ChS]}  J. Cheeger and J. Simons, {\sl
Differential Characters and Geometric Invariants},
in Geometry and Topology, Lecture Notes in Math. {\bf 1167},
Springer--Verlag, New York, 1985, 50--80.

\item {[CJS] }   R.L. Cohen, J. D. S. Jones and G. B. Segal,   {\sl
Morse theory and classifying spaces }, Stanford Preprint, 1993.

\item {[deR]} G. de Rham,   Vari\'et\'es Diff\'erentiables,
Hermann, Paris, 1973.

\item {[ES] } G. Ellingsrud and S. A. Stromme, {\sl  Toward the Chow
ring of the Hilbert scheme of $\bbp^2$, } J. Reine Agnew. Math. {\bf
441 } (1993 ),  33-44.

\item {[F]} H. Federer,   Geometric Measure Theory,
Springer Verlag, New York, 1969.

\item {[Fr]}  T. Frankel, {\sl Fixed points and torsion on  \Kahler
manifolds},
  Ann.of Math.,{\bf 70} (1959), 1-8.

\item {[GS] }  H. Gillet and C. Soul\'e,
{\sl Characteristic classes for algebraic vector bundles with
Hermitian metrics I, II}, Ann. of Math., {\bf 131 } (1990 ), 163-203,
205-238.

\item{[HL$_1$]}  F.R. Harvey and H.B. Lawson, Jr.,  A Theory of
 Characteristic Currents Associated with a Singular Connection,
   Ast\'erisque {\bf 213},  Soci\'et\'e   Math. de France, Paris, 1993.

\item {[HL$_2$]} --- , {\sl Singularities and Chern-Weil theory, I --
The local MacPherson formula}, Asian J. Math., {\bf 4}, No. 1
(2000), 71--96.

\item {[HL$_3$]} --- , {\sl Singularities and Chern-Weil theory, II --
Geometric atomicity}, (to appear).

\item {[HL$_4$]} --- , {\sl Morse Theory and Stokes' Theorem},
in  Surveys in Differential Geometry, VII,  International Press, to appear.

\item {[HL$_5$]} --- , {\sl Poincar\'e-Pontrjagin duality for
differential characters}, (to appear).

\item {[HP]}  R. Harvey and J. Polking, {\sl Fundamental solutions in
complex analysis, Part I}, Duke Math. J. {\bf 46} (1979), 253-300.

\item {[K]}  M. Kontsevich, {\sl Feynman diagrams and
low-dimensional topology},  First European Congress of
Mathematics, Vol. II, Birh\"auser Verlag, Basel, 1994,  pp. 97-122.

\item {[L]}  J. Latschev,  Gradient flows of Morse-Bott functions,
{\it Math. Ann.\/} {\bf 318}  (2000), 731--759.

\item {[La]} F. Laudenbach, {\sl On the Thom-Smale complex}, in
An Extension of a Theorem of  Cheeger and M\"uller,
by J.-M. Bismut and W. Zhang, Ast\'erisque, v. 205:  S.M.F., 1992.

\item {[Mac$_1$]} R. MacPherson, {\sl  Singularities of vector bundle
 maps\/},  {\it Proc.\  of Liverpool  Singularities Symposium, I},
Springer Lecture Notes in Mathematics, {\bf 192} (1971), 316-318.

\item {[Mac$_2$]} --- , {\sl  Generic vector bundle maps } in
`` Dynamical Systems, Proceedings of Symposium -- University of Bahia,
Salavador 1971'',  Academic Press, New York, 1973, pp. 165-175.

\item {[S]} S. Smale, {\sl On gradient dynamical systems} , Annals of
Math. {\bf 74} (1961), 199-206.

\item {[So]}  A. J. Sommese, {\sl Extension theorems for reductive
groups
 actions on compact \Kahler manifolds},  Math. Ann., {\bf 218}
(1975), 107-116.

\item {[T]}  R. Thom, {\sl Sur une partition en cellules associ\'ees
\`a une fonction sur une vari\'et\'e}, CRAS  {\bf 228} (1949), 973-975.

\item {[W]} E. Witten, {\sl Supersymmetry and Morse theory  },
  J. Diff. Geom., {\bf 17  }, no. 4 (1982 ),  661-692.

\end